\documentclass[12pt]{article}

\usepackage{amssymb,amsmath}

\usepackage[francais,english]{babel}
\usepackage[T1]{fontenc}
 \makeatletter \@addtoreset{equation}{section}
\newtheorem{thm}{Theorem}[section]
\newtheorem{lem}{Lemma}[section]
\newtheorem{cor}{Corollary}[section]
\newtheorem{claim}{Proposition}[section]

\begin{document}

\title{Hopf maximum
principle violation for elliptic equations with non-Lipschitz
nonlinearity}

\author{Y. Il'yasov \footnote{The first author was 
partly supported by grants  RFBR 08-01-00441-a,
08-01-97020-p-a}\\Institute of Mathematics RAS,\\
Ufa, Russia\and Y. Egorov\\
Universit\'e Paul Sabatier,\\
Toulouse, France}
\date{$\ $}
\maketitle

\abstract{ We consider elliptic equations with non-Lipschitz
nonlinearity $$ -\Delta u = \lambda |u|^{\beta-1}u-|u|^{\alpha-1}u$$ in a smooth
bounded domain $\Omega \subset \mathbb{R}^n$, $n\geq 3$, with
Dirichlet boundary conditions; here $0<\alpha<\beta<1$. We prove the existence of a weak nonnegative solution
which does not satisfy  the Hopf boundary maximum principle,
provided that $\lambda$ is large enough and
$n>2(1+\alpha) (1+\beta)/(1-\alpha)(1-\beta)$.}

\section{Introduction}

Let $\Omega$ be a bounded domain in $\mathbb{R}^n$, $n\geq 3$ with
a smooth boundary $\partial \Omega$, which is strictly star-shaped
with respect to the origin in $\mathbb{R}^n$. We consider the
following problem:
\begin{equation}
\label{10R} \left\{
\begin{array}{l}
\ -\Delta u = \lambda |u|^{\beta-1}u-|u|^{\alpha-1}u~~ \mbox{in}~ \Omega,
 \\ \\
\ u =0~ \mbox{on}~ \partial\Omega.
%\label{1.2}
%
\end{array}
\right.
\end{equation}
Here $\lambda$ is a real parameter and $
0<\alpha<\beta<1$, so that the nonlinearity $f(\lambda,u):=\lambda |u|^{\beta-1}u-|u|^{\alpha-1}u$  on the right-hand side of (\ref{10R}) is non-Lipschitzean at zero.

Our interest to this problem has been induced by investigations of
J.I. D\'iaz, J. Hern\'andez in paper \cite{diaz}. In case of
dimension $n=1$  when  $\Omega=(-1,1)$, among other results, they
showed that  for certain values $\lambda>0$ equation (\ref{10R})
possesses solutions $u(x)$, $x\in (-1,1)$ with a special feature
\begin{equation}\label{dher}
u(-1)=u(1)=0,~~u'(-1)=u'(1)=0.
\end{equation}
This means a Hopf boundary maximum principle violation on $x=-1$,
$x=1$ and a loss of the uniqueness for initial value problem to
(\ref{10R}) with  $u(-1)=u'(-1)=0$, since $u\equiv 0$ satisfies
also to (\ref{10R}).  Furthermore, it can easily be shown that the
existence of such a solution with $\lambda_0>0$ yields the
existence of a  set continuum nonnegative solutions of this
boundary value problem for any $\lambda>\lambda_0$. Observe that
property (\ref{dher}) implies that a function $u$  is also a
weakly solution of (\ref{10R}) on the whole line $\mathbb{R}$.
Note that when the  nonlinearity $f(\lambda,u)$ is a locally
Lipschitz function such a phenomenon is impossible due to the
uniqueness solution of initial value problem and/or a Hopf
boundary maximum principle.
 
This rise a question {\it as to whether the similar phenomena may be
occurred in case of the higher dimensions $n>1$}. More precisely
whether the Hopf boundary maximum principle  holds for
(\ref{10R}) when $n>1$ and  the nonlinearity
$f(\lambda,u)$ is non-Lipschitz. To find an answer to this
question is a main goal in the present work.

Let us state our main result.
We consider a weak solution $u \in H^1_0:=H^1_0(\Omega)$, where $H^1_0(\Omega)$ denotes the closure $C^\infty_0(\Omega)$ in standard Sobolev space $H^1(\Omega)$ with the norm $||\cdot||_1$. We say that a weak solution $u \in H^1_0$ of (\ref{10R}) is a non-regular, if $u \in C^1(\overline{\Omega})$ and $\frac{\partial u}{\partial \nu} =0$ on $\partial \Omega$.

Our main result is the following
\begin{thm}\label{Th1}
Let $\Omega$ be a bounded domain in $\mathbb{R}^n$,
$n\geq 3$ with smooth boundary, which is strictly star-shaped with respect to the
origin. Assume that $0<\alpha<\beta<1$ and
$n>2(1+\alpha)(1+\beta)/(1-\alpha)(1-\beta)$. Then there exists
$\lambda^*>0$ such that for all $\lambda \geq \lambda^*$ problem
(\ref{10R}) has  a non-regular solution solution $u_\lambda$, which is nonnegative in $\Omega$. Moreover, the number of such solutions for $\lambda > \lambda^*$ is infinite.

Furthermore, for any $\lambda > \lambda^*$ problem (\ref{10R}) has a weak solution $w_\lambda \in C^1(\overline{\Omega})$,  which is nonnegative in $\Omega$ but is not non-regular solution, i.e. $\frac{\partial
w_\lambda}{\partial \nu} \neq 0$ on some subset $U \subseteq \partial \Omega$
of positive $(n-1)$-dimensional Lebesgue measure.
\end{thm}
The proof of the theorem relies on the variational arguments.
Furthermore, basic ingredients in the proof consist in using  Pohozaev's
identity \cite{poh} corresponding to (\ref{10R}) and in applying the spectral analysis with respect to the fibering procedure introduced in \cite{ilyas}.

{\bf Remark 1.1.}\, If one considers the radial symmetric solutions of (\ref{10R}) in 
the ball $B_R$, then the second part of Theorem \ref{Th1} implies that the  weak radial solution $w_\lambda$ of (\ref{10R}) is positive in the ball $B_R$ and satisfies  the Hoph boundary maximum principle, i.e. $\frac{\partial u_\lambda}{\partial \nu}(R) <0$.

{\bf Remark 1.2.}\, In the theory of integrable systems, the non-regular type solutions are known as the compactons: solitary waves with compact support \cite{ros}.

The paper is organized as follows. In Section 2,  we apply the spectral analysis related with the
fibering procedure \cite{ilyas} to introduce two spectral points $\Lambda_0,\Lambda_1$
which play basic role in the proof of the main result. In Section 3, we derive some important consequences from Pohozaev's identity.  In Section 4, we prove  existence of the solution to an auxiliary constrained minimization problem. In Section 5, we prove Theorem \ref{Th1}. Section 6 is an appendix where  some technical result is
proved.

\section{Spectral analysis with respect to the fibering procedure}
In this section we apply the spectral analysis with respect to the
fibering procedure \cite{ilyas} to introduce two spectral points
which will play important roles in the proof of the main result.

Observe that problem (\ref{10R}) is the Euler-Lagrange equation of
the functional
\begin{eqnarray}
E_\lambda(u)=\frac{1}{2}T(u) -\lambda\frac{1}{\beta+1}B(u) +
\frac{1}{\alpha+1}A(u),    \label{euler}
\end{eqnarray}
where we use the notations
\begin{eqnarray*}
T(u)=\int_\Omega |\nabla u|^2\,dx,~ B(u)=\int_\Omega
u^{\beta+1}\,dx,~ A(u)=\int_\Omega u^{\alpha+1}\,dx.
\end{eqnarray*}

Let $u\in H^1_0$. Consider the function
$e_\lambda(t):=E_\lambda(tu)$ defined for $t\in \mathbb{R}^+$.
Introduce the functionals $
Q_\lambda(u)=e_\lambda'(t)|_{t=1},~~L_\lambda(u)=e_\lambda''(t)|_{t=1}
$ for $u \in H^1_0$. Then
$$
Q_\lambda(u)=\ T(u) -\lambda B(u) + A(u),~L_\lambda(u):=T(u)
-\lambda \beta B(u) + \alpha A(u).
$$

Let $u \in H^1_0\setminus \{0\}$. Following the spectral analysis
 \cite{ilyas}, we solve the system
\begin{equation}
\label{sis0}  \left\{
\begin{array}{l}
Q_\lambda(tu)=t^2T(u) -\lambda t^{1+\beta} B(u)
+ t^{1+\alpha}A(u)=0 \\ \\
E_\lambda(tu)=\frac{t^2}{2}T(u) -\lambda
\frac{t^{1+\beta}}{1+\beta}B(u) + \frac{t^{1+\alpha}}{1+\alpha}
A(u)=0
\end{array}
\right.
\end{equation}
and find the corresponding solution
\begin{equation}\label{P1}
  t_0(u)=\left(\frac{2(\beta-\alpha)}{(1+\alpha)(1-\beta)}
  \frac{A(u)}{T(u)}\right)^\frac{1}{1-\alpha},
\end{equation}
$$
\lambda_0(u)=c^{\alpha,\beta}_0 \lambda(u),
$$
where
$$
 c^{\alpha,\beta}_0=\frac{(1-\alpha)(1+\beta)}{(1-\beta)(1+\alpha)}
\left(\frac{(1+\alpha)(1-\beta)}{2(\beta-\alpha)}\right)^\frac{\beta-\alpha}{1-\alpha}
$$
and
\begin{equation}\label{lam}
\lambda(u)=
\frac{A(u)^\frac{1-\beta}{1-\alpha}T(u)^\frac{\beta-\alpha}{1-\alpha}}{B(u)}.
\end{equation}

Thus with respect to the fibering procedure, we have the following
spectral point
\begin{equation}\label{P3}
\Lambda_0=\inf_{H^1_0\setminus \{0\}}\lambda_0(u).
\end{equation}

Introduce the second  point $\Lambda_1$.  Let $u \in
H^1_0\setminus \{0\}$. Consider now the following system
\begin{equation}
\label{sis2} \left\{
\begin{array}{l}
Q_\lambda(tu)=t^2T(u) -\lambda t^{1+\beta} B(u)
+ t^{1+\alpha}A(u)=0 \\ \\
L_\lambda(tu)=t^2T(u) -\lambda \beta t^{1+\beta}B(u) + \alpha
t^{1+\alpha}A(u)=0
\end{array}
\right.
\end{equation}
for $t\in \mathbb{R}^+$, $\lambda\in \mathbb{R}^+$. Solving this
system  we find as above
\begin{equation}\label{P20}
\lambda_1(u)=c^{\alpha,\beta}_1 \lambda(u),
\end{equation}
where
\begin{equation}\label{P40}
 c^{\alpha,\beta}_1=\frac{1-\alpha}{1-\beta}
\left(\frac{1-\beta}{\beta-\alpha}\right)^\frac{\beta-\alpha}{1-\alpha}.
\end{equation}
Then we have
\begin{equation}\label{P30}
\Lambda_1=\inf_{H^1_0\setminus \{0\}}\lambda_1(u).
\end{equation}
\begin{claim}\label{proL}
$0<\Lambda_1<\Lambda_0<+\infty$.
\end{claim}
{\it Proof.}\, Observe that $\lambda_1(u)=
C^{\alpha,\beta}\lambda_0(u)$ for any $u\in H^1_0\setminus \{0\}$
with $C^{\alpha,\beta}= c^{\alpha,\beta}_0/ c^{\alpha,\beta}_1$.
It is not hard to show that $C^{\alpha,\beta}<1$, therefore
$\Lambda_1<\Lambda_0$.

It is clear that $\Lambda_0<+\infty$. Let us show that
$0<\Lambda_1$. Note that $\lambda(u)$ is a zero-homogeneous
function on $H^1_0\setminus \{0\}$. Therefore we may restrict the
infimum in (\ref{P30}) to the set $S:=\{v \in H^1_0:
||v||_1=1\}$.

Set
\begin{equation}\label{P5}
    \gamma=\frac{(1+\alpha)(2^*-1-\beta)}{(2^*-1-\alpha)},~
    p=\frac{(1+\alpha)}{\gamma},~q=\frac{2^*}{(1+\beta-\gamma)},
\end{equation}
where $2^*=2n/(n-2)$. Then  $p,q>1$, $1/p+1/q=1$ and by the H\"{o}lder
inequality we have
\begin{equation}\label{P6}
 B(u)\leq (\int_\Omega u^{2^*}dx)^{1/q}\cdot A(u)^{1/p}.
\end{equation}
By Sobolev's inequality
$$
\int_\Omega u^{2^*}dx\leq C_0||u||_1^{2^*}=C_0<+\infty.
$$
for $u \in S$, where $C_0$ does not depend on $u \in H^1_0$. Hence
for any $u\in S$ we have
\begin{equation}
\lambda(u)=\frac{A(u)^{(1-\beta)}}{B(u)^{(1-\alpha)}} \geq
c_{0}A(u)^{-\frac{(2^*-2)(\beta-\alpha)}{(2^*-1-\alpha)}},
\label{PPP}
\end{equation}
where $0<c_0<+\infty$ does not depend on $u \in H^1_0$. Since
$A(u)\leq C_1<+\infty$ on $S$ we see from (\ref{P3}) that
$\Lambda_1>0$.
\hspace*{\fill}\rule{3mm}{3mm}\\
\par\noindent

\section{Pohozaev's identity}

We will need the following regularity result
\begin{claim}\label{Reg}
Assume that $0<\alpha<\beta<1$. Suppose that $u \in H^1_0$ is a weak
solution of (\ref{10R}). Then $u \in
C^{1,\kappa}(\overline{\Omega})$ for $\kappa \in (0,1)$.
\end{claim}
{\it Proof.} Let $u \in H^1_0$ be a weak solution of (\ref{10R}).
Since $|f(\lambda,u)|<C(1+|u|)$, $u \in
\mathbb{R}$ with some $C>0$, then (see Lemma B.3 in \cite{str}) $u
\in L^q(\Omega)$ for any $q<\infty$. This implies that  $-\Delta u =
f(\lambda,u) \in L^q(\Omega)$ for any $q<\infty$.
Thus, by the Cald\'eron-Zygmund inequality (see \cite{giltrud})
$u\in H^{2,q}(\Omega)$, whence $u \in
C^{1,\kappa}(\overline{\Omega})$ for $\kappa \in (0,1)$ by the
Sobolev embedding theorem.
\hspace*{\fill}\rule{3mm}{3mm}\\
\par\noindent

We will denote by $P_\lambda$  the functional
$$
P_\lambda(u):=\displaystyle{\frac{(n-2)}{2n}T(u)
-\lambda\frac{1}{\beta+1}B(u) + \frac{1}{\alpha+1}A(u)}
$$
defined for $u \in H^1_0$.

\begin{lem}\label{lem1}
Suppose that $\Omega$ is a smooth bounded domain in $\mathbb{R}^n$,
$n\geq 3$, which is strictly star-shaped with respect to the
origin in $\mathbb{R}^n$. Let $u$ be a weak solution of
(\ref{10R}), $u \in H^1_0$. Then the following Pohozaev identity
holds
\begin{equation}\label{poh}
P_\lambda(u)+ \frac{1}{2n}\int \left|\frac{\partial u}{\partial
\nu}\right|^2\,x\cdot \nu\,dx=0.
\end{equation}
\end{lem}
{\it Proof.}\, By Proposition \ref{Reg} we know that $u \in
H^{2,2}\cap C^1\cap H^1_0$. Thus, since $f(\lambda,u)$ is a continuous function on $\mathbb{R}$,
 we are in position to apply Lemma 1.4  in \cite{str}, that completes 
 the proof.
\hspace*{\fill}\rule{3mm}{3mm}\\
\par\noindent

Let $u \in H^1_0$. Based on the ideas of the  spectral analysis
with respect to the fibering procedure \cite{ilyas} we consider
the following system of equations
\begin{equation}
\label{sis1} \left\{
\begin{array}{l}
Q_\lambda(u):=\ T(u) -\lambda B(u)
+ A(u)=0 \\ \\
L_\lambda(u):=T(u) -\lambda \beta B(u) + \alpha A(u)=0
\\ \\
P_\lambda(u):=\displaystyle{\frac{(n-2)}{2n}T(u)
-\lambda\frac{1}{\beta+1}B(u) + \frac{1}{\alpha+1}A(u)=0}.
%\label{1.2}
%
\end{array}
\right.
\end{equation}
The computation of the corresponding determinant shows that
this system is  solvable if and only if
\begin{equation}\label{theta}
\theta\equiv 2(1+\alpha)(1+\beta)-n(1-\alpha)(1-\beta)=0.
\end{equation}
Note that $\theta<0$ if and only if 
$$
n>2(1+\alpha)(1+\beta)/(1-\alpha)(1-\beta).
$$
Observe that the equation $e_\lambda'(t)=0$, $t>0$ has at most two roots
$t^1(u):=t^1_\lambda(u)\in \mathbb{R}^+$ and $t^2(u):=t_\lambda^2(u)\in \mathbb{R}^+$ such that $t^1(u)\leq t^2(u)$, $e_\lambda''( t^1(u))\leq 0$ and $e_\lambda''(t^2(u))\geq 0$.

\begin{claim}\label{pro}
Assume that $\theta <0$. If $u \in H^1_0\setminus \{0\}$ and $t>0$ are
such that $Q_\lambda(tu)=0$ and $P_\lambda(tu)\leq 0$ then we have
$$
 L_\lambda(tu) >0.
$$
\end{claim}
{\it Proof.}\, Let $u \in H^1_0\setminus \{0\}$ and $t>0$ as in the assumption.
Then
\begin{eqnarray}
&&T(u)=\lambda t^{\beta-1}B(u)-t^{\alpha-1}A(u)\label{ver1}\\
&&L_\lambda(tu) >0~\Leftrightarrow~\lambda t^{\beta-\alpha}\frac{(1-\beta)}{(1-\alpha)}B(u)> A(u).\label{ver2}
\end{eqnarray}
Equality (\ref{ver1}) implies that  $P_\lambda(t^1(u)u)\leq 0$
holds if and only if
\begin{equation}\label{ver3}
\lambda t^{\beta-\alpha}\frac{[2(1+\beta)+n(1-\beta)](1+\alpha)}
{[2(1+\alpha)+n(1-\alpha)](1+\beta)}B(u)\geq A(u).
\end{equation}
Observe, that the inequality $\theta <0$ implies
\begin{equation}\label{ver4}
\frac{[2(1+\beta)+n(1-\beta)](1+\alpha)}
{[2(1+\alpha)+n(1-\alpha)](1+\beta)}<  \frac{(1-\beta)}{(1-\alpha)}.
\end{equation}
Thus (\ref{ver4}) and (\ref{ver3}) give
$$
\lambda t^{\beta-\alpha}\frac{(1-\beta)}{(1-\alpha)}B(u)> A(u)
$$
and therefore by (\ref{ver2})  the proof is complete.
\hspace*{\fill}\rule{3mm}{3mm}\\
\par\noindent

\begin{cor}\label{CF}
If $u_0$ is a non-regular solution solution of (\ref{10R}) then
$E_\lambda(u_0)>0$. Furthermore, if in addition $\theta<0$, then
$$
Q_\lambda(u_0)=0,~P_\lambda(u_0)=0,~L_\lambda(u_0)>0.
$$
\end{cor}
{\it Proof.}\, Observe that if $u_0$ is the non-regular solution solution of
(\ref{10R}), then by (\ref{poh}) we have $ P_\lambda(u_0)=0$.
Hence using  $E_\lambda(u)=P_\lambda(u)+(1/2n)T(u)$ we get
$$
E_\lambda(u_0)=\frac{1}{n}T(u_0)>0.
$$
Note that $Q_\lambda(u_0)=0$ if $u_0$ is a solution of (\ref{10R}),
and  $P_\lambda(u_0)=0$ if in addition this solution is the non-regular solution solution.
Hence assumption $\theta<0$ and Proposition \ref{pro} imply that
$L_\lambda(u_0)>0$.

\hspace*{\fill}\rule{3mm}{3mm}\\
\par\noindent

\section{Constrained minimization problems}
Consider the following constrained minimization problem:
\begin{equation}
\label{min1} \left\{
\begin{array}{l}
\ E_\lambda(u) \to \min \\ \\
Q_\lambda(u)= 0.
%\label{1.2}
%
\end{array}
\right.
\end{equation}
We denote by
$$
M_\lambda:=\{w \in H^1_0: ~Q_\lambda(u)=0 \}
$$
the admissible set of (\ref{min1}), and by
$\hat{E_\lambda}:=\min\{E_\lambda(u):~ u \in M_\lambda\}$ the minimal
value in this
 problem. We say that $(u_m)$ is a
minimizing sequence of (\ref{min1}), if
\begin{equation}\label{maxs}
E_\lambda(u_m) \to \hat{E_\lambda}~~\mbox{as}~~m\to
\infty~~\mbox{and}~ u_m \in M_\lambda,~m=1,2,...
\end{equation}

\begin{claim}\label{pro4}
If  $\lambda>\Lambda_1$, then the set $M_\lambda$ is not empty,
meanwhile the set $M_\lambda$ is empty when $\lambda<\Lambda_1$.
\end{claim}
{\it Proof.}\, Let $\lambda>\Lambda_1$. Then by (\ref{P3}) there
exists $u \in H^1_0\setminus \{0\}$ such that
$\Lambda_1<\lambda(u)<\lambda$ and $L_{\lambda(u)}(t(u)u)=
0,~Q_{\lambda(u)}(t(u)u)= 0$. Hence, $Q_\lambda(t(u)u)< 0$, since
$\lambda(u)<\lambda$ and therefore there exists $t>0$ such that
$Q_\lambda(tu)= 0$, i.e. $tu \in M_\lambda$.

The proof of the second part of the Proposition follows
immediately from the definition (\ref{P3}) of $\Lambda_1$.

\hspace*{\fill}\rule{3mm}{3mm}\\
\par\noindent

From here it follows that
\begin{cor}
$\hat{E_\lambda}<+\infty$ for any $\lambda> \Lambda_1$.
\end{cor}

\par\noindent

\subsection{Existence of the solution of (\ref{min1}).}

\begin{lem}\label{le1e}
For any $\lambda> \Lambda_1$ problem (\ref{min1}) has  a solution
$u_0 \in H^1_0\setminus \{0\}$, i.e.
$E_\lambda(u_0)=\hat{E_\lambda}$ and $u_0\in
M_\lambda$.
\end{lem}
{\it Proof.}\, Let $\lambda>\Lambda_1$. Then $M_\lambda$ is not empty  and there is  a
minimizing sequence $(u_m)$ of (\ref{min1}). Set $t_m\geq 0$ and $v_m\in
H^1_0$, $m=1,2,...,$ such that $u_m=t_m v_m$, $||v_m||_1=1$.

Let us show that $\{t_m\}$ is bounded.
Observe that
\begin{eqnarray}\label{Eg1}
1 -\lambda t^{\beta-1}_m B(v_m)
+ t^{\alpha-1}_m A(v_m)=0,
\end{eqnarray}
since $Q_\lambda(t_mu_m)=0$, $m=1,2,...$. Note that since $||v_m||_1=1$, $B(v_m),
A(v_m)$ are bounded.

Suppose  that there exists a subsequence again denoted $(t_m)$ such that   $t_m \to \infty$ as $m\to
+\infty$. Then the left hand side of
(\ref{Eg1}) tends to 1 as $m\to +\infty$  what contradicts to
the assumption $Q_\lambda(u_m)=0$, $m=1,2,...$.

Suppose now  that there exist subsequences again denoted $(t_m)$, $(v_m)$ such that $t_m \to 0$ and/or $v_m \to
0$ weakly in $H^1_0$ as $m\to +\infty$.

Assume that $t^{\alpha-1}_m A(v_m) \to C$ as $m\to \infty$, where $0\leq C<+\infty$.
Then $\lambda t^{\beta-1}_m B(v_m) \to 1+C$ as $m\to \infty$. By (\ref{P6})
we have
$ B(v_m)\leq C_0\cdot A(v_m)^{1/p}$,
where $0<C_0<+\infty$ does not depend on $m=1,2,...$. Therefore
\begin{equation}\label{P7}
 t_m^{\beta-1}B(v_m)\leq C_0\cdot t_m^{\beta-1}A(v_m)^{1/p}=t_m^{\beta-1+\frac{(1-\alpha)}{p}}(t_m^{\alpha-1}A(v_m))^{1/p}.
\end{equation}
Let us show that
\begin{equation}\label{Eg3}
\beta-1+\frac{(1-\alpha)}{p}>0.
\end{equation}
Substituting $p=\frac{(2^*-1-\alpha)}{(2^*-1-\beta)}$ we get
\begin{eqnarray*}
&&\beta-1+\frac{(1-\alpha)}{p}=\frac{(\beta-1)(2^*-1-\alpha)+(1-\alpha)(2^*-1-\beta)}{(2^*-1-\alpha)}.
\end{eqnarray*}
Since
\begin{eqnarray*}
&&(\beta-1)(2^*-1-\alpha)+(1-\alpha)(2^*-1-\beta)=\\
&&2^*(\beta-\alpha)-(1+\alpha)(\beta - 1)-(1-\alpha)(1+\beta)=\\
&&2^*(\beta-\alpha)-2(\beta-\alpha)=(\beta-\alpha)(2^*-2)>0,
\end{eqnarray*}
we get the desired conclusion.  Hence  the right hand side in (\ref{P7}) tends to zero and therefore
$t^{\beta-1}_m B(v_m) \to 0$ as $m\to \infty$, which contradicts our assumption.

Assume now that $t^{\alpha-1}_m A(v_m) \to +\infty$ as $m\to \infty$.
Then by (\ref{Eg1}) we have
\begin{equation}\label{Eg2}
\frac{A(v_m)}{\lambda t^{\beta-\alpha}_m B(v_m)} \to  1
\end{equation}
as $m\to \infty$.
Using (\ref{P6}) we deduce
\begin{eqnarray}\label{Eg4}
\frac{A(v_m)}{\lambda t^{\beta-\alpha}_m B(v_m)}>c_0\frac{A(v_m)^{(p-1)/p}}{ t^{\beta-\alpha}_m }=c_0\frac{(t_m^{(\alpha-1)}A(v_m))^{(p-1)/p}}{t_m^{(p-1)(\alpha-1)/p+\beta-\alpha}},
\end{eqnarray}
where $0<c_0<+\infty$ does not depend on $m=1,2,...$.
Using (\ref{Eg3}) we get
$$
\frac{(p-1)(\alpha-1)}{p}+\beta-\alpha=\beta-1+\frac{1-\alpha}{p}>0.
$$
This implies that the right hand side of (\ref{Eg4}) tends to $+\infty$, contrary to (\ref{Eg2}).

Thus  $(u_m)$ is
bounded in $H^1_0$, and hence by Sobolev's embedding theorem, $(u_m)$
has a subsequence which  converges weakly in $H^1_0$ and strongly in $L_p$,
$1<p<2^*$ . Denoting this subsequence again by
$(u_m)$ we get $u_m \to u_0$ weakly in $H^1_0$ and strongly
in $L_p$, $1< p<2^*$ for some $u_0\in H^1_0$. By the
above, the sequences $(t_m)$ and $(v_m)$ are separated from zero
and therefore $u_0 \neq 0$. Thus $E_\lambda(u_0)\leq
\hat{E}_\lambda$ and $Q_\lambda(u_0) \leq 0$. Assume
$Q_\lambda(u_0) < 0$. Then
$Q_\lambda(t^2_\lambda(u_0)u_0) = 0$, i.e.
$t^2_\lambda(u_0)u_0\in M_\lambda$ and
$E_\lambda(t^2_\lambda(u_0)u_0)<E_\lambda(u_0)\leq
\hat{E}_\lambda$. Hence we get a contradiction and therefore
$E_\lambda(u_0)= \hat{E}_\lambda$ and $Q_\lambda(u_0)
= 0$.  This completes the proof of Lemma \ref{le1e}.

\hspace*{\fill}\rule{3mm}{3mm}\\
\par\noindent

From the definition (\ref{P30}) of $\Lambda_0$ and using arguments as in the proof of Lemma
\ref{le1e} it is not hard to derive
\begin{cor}\label{proE}
If  $\lambda>\Lambda_0$, then $\hat{E_\lambda}<0$.
If  $\Lambda_1<\lambda<\Lambda_0$, then $0<\hat{E_\lambda}<+\infty$, and if $\lambda=\Lambda_0$, then $\hat{E_\lambda}=0$.
\end{cor}

\subsection{Existence of the solution of (\ref{10R}).}

Let $\lambda>\Lambda_1$ then by Lemma \ref{le1e} there exists a
solution $u_0 \in H^1_0\setminus \{0\}$ of (\ref{min1}). This implies that there
exist Lagrange multipliers $\mu_1$, $\mu_2$ such that
\begin{equation}\label{eq2}
\mu_1DE_\lambda(u_0)=\mu_2 DQ_\lambda(u_0),
\end{equation}
 and $|\mu_1|+|\mu_2|\neq 0$.
\begin{claim}\label{Lag}
Let $\theta<0$, $\lambda>\Lambda_1$ and
 $u_0 \in H^1_0$ be a solution of (\ref{min1}).
Assume that $P_\lambda(u_0)\leq 0$. Then  $u_0$ is a weak
nonnegative solution of (\ref{10R}).
\end{claim}
{\it Proof.}\, Note that by Proposition \ref{pro} we have $L_\lambda(u_0)\neq 0$, since $\theta<0$, $Q_\lambda(u_0)=0$ and by the assumption $P_\lambda(u_0)\leq 0$. From (\ref{min1}) and (\ref{eq2}) we have
$0=\mu_1Q_\lambda(u_0)=\mu_2 L_\lambda(u_0)$. But $L_\lambda(u_0)\neq 0$ and therefore
$\mu_2=0$. Thus by (\ref{eq2}) we have $DE_\lambda(u_0)=0$. Since
$E_\lambda(|u_0|)=E_\lambda(u_0)$,
$Q_\lambda(|u_0|)=Q_\lambda(u_0)=0$ we may assume that
$u_0\geq 0$.   This completes the proof.
\hspace*{\fill}\rule{3mm}{3mm}\\
\par\noindent

\begin{cor}\label{Lag2}
Let  $\theta<0$, $\lambda\geq\Lambda_0$ and
 $u_0 \in H^1_0$ be a solution of (\ref{min1}).
Then $P_\lambda(u_0)<0$ and $u_0$ is a weak solution of
(\ref{10R}).
\end{cor}
{\it Proof.}\, Corollary \ref{proE} implies  $\hat{E}_\lambda <0$
when $\lambda>\Lambda_0$, and $\hat{E}_\lambda =0$ when
$\lambda=\Lambda_0$. Hence for any $\lambda\geq\Lambda_0$ we have
$E_\lambda(u_0)\leq 0$ and therefore the identity
$E_\lambda(u_0)=P_\lambda(u_0)+(1/2n)T(u_0)$ implies that
$P_\lambda(u_0)<0$. Applying now Proposition \ref{Lag} we complete
the proof.
\hspace*{\fill}\rule{3mm}{3mm}\\
\par\noindent

\section{Proof of Theorem \ref{Th1}}
Let us introduce
\begin{equation}\label{Pr1}
Z:=\{\lambda>0:~\exists u_\lambda \in
M_\lambda~\mbox{s.t.}~E_\lambda(u_\lambda)=\hat{E}_\lambda,~P_\lambda(u_\lambda)<0\}.
\end{equation}
By assumption,  $n>2(1+\alpha)(1+\beta)/(1-\alpha)(1-\beta)$, i.e. $\theta<0$. 
Hence Lemma \ref{le1e} and Corollary \ref{Lag2} imply that $Z$ is
bounded below by $\Lambda_1$ and $[\Lambda_0, +\infty)\subset Z$,
i.e. $Z\neq \emptyset$. Furthermore,  Lemma \ref{app} from
Appendix yields that the maps $G_{(\cdot)}(u_{(\cdot)}): \lambda
\mapsto P_\lambda(u_\lambda)$, $E_{(\cdot)}(u_{(\cdot)}): \lambda
\mapsto E_\lambda(u_\lambda)$ are continuous functions in
$(\Lambda_1, +\infty)$ and hence $Z\cap (\Lambda_1, +\infty)$ is
an open set in $\mathbb{R}$.

Introduce
$$
\lambda^*:=\inf Z.
$$
\begin{lem}\label{Pr2}
There exists a solution $u^*$ of (\ref{min1}) with
$\lambda=\lambda^*$. Furthermore, $\Lambda_1<\lambda^*$ and
$P_{\lambda^*}(u^*)=0$.
\end{lem}
{\it Proof.}\, Since $Z$ is an open set, we can find a sequence
$\lambda_m \in Z$, $m=1,2,...$ such that $\lambda_m \to
\lambda^*$ as $m \to \infty$. By definition of $Z$ for any
$m=1,2,...$ there exists solution $u_{\lambda_m}$ of (\ref{min1})
such that $P_{\lambda_m}(u_{\lambda_m})<0$.  Lemma \ref{app} from
Appendix yields the existence of the limit solution $u^*$ of
(\ref{min1}) and the existence of a subsequence (again denoted by $(u_{\lambda_m})$) such that $u_{\lambda_m} \to u^*$ strongly
 in $H^1$ as $\lambda_m \to \lambda^*$. This yields by continuity that $P_{\lambda^*}(u^*)\leq 0$.

Let us show that $\Lambda_1<\lambda^*$. To obtain a contradiction suppose, that $\Lambda_1=\lambda^*$. Then by the proof of Lemma \ref{app} from
Appendix A we know that $\Lambda_1=\lambda_1(u^*)$. Thus $u^*$ is a critical point of $\lambda(u)$. This implies that $t_0(u^*)u^*$ is a weak solution of (\ref{10R}) with $\lambda=\Lambda_0$. Note that since
$P_{\lambda_m}(u_{\lambda_m})<0$, $m=1,2,...$,  by Proposition \ref{Lag} $u_{\lambda_m}$ weakly satisfies (\ref{10R}) with $\lambda=\lambda_m$, $m=1,2,...$.
From here and since $u_{\lambda_m} \to u^*$ strongly
 in $H^1$ as $\lambda_m \to \Lambda_1$, we derive that $u^*$ is a weak solution of (\ref{10R}) with $\lambda=\Lambda_1$. But $\Lambda_1\neq \Lambda_0$ and we get a contradiction.

Thus $\lambda^*\in (\Lambda_1, +\infty)$ and  $Z$ is an open subset in $(\Lambda_1, +\infty)$. Suppose, contrary to our claim, that  $P_{\lambda^*}(u^*)< 0$. Then  $\lambda^* \in Z$. However, since $Z$ is an open set, this is impossible, and therefore $P_{\lambda^*}(u^*)= 0$.

\hspace*{\fill}\rule{3mm}{3mm}\\
\par\noindent
{\it Conclude of the proof of Theorem \ref{Th1}}
\par\noindent

By Proposition \ref{Lag}, $u^*$ is a weak nonnegative solution of
(\ref{10R}) and by Proposition \ref{Reg} $u^* \in C^{1,\kappa}(\overline{\Omega})$ for
$\kappa \in (0,1)$.  Hence by Lemma \ref{lem1} Pohozaev's identity
holds and whence $\int \left|\frac{\partial u^*}{\partial
\nu}\right|^2\,x\cdot \nu\,dx=0$, since $P_{\lambda}(u^*)= 0$. By the assumption $\Omega$ is a
strictly star-shaped domain with respect to the origin, i.e.
$x\cdot \nu>0$ on $\partial \Omega$. Therefore $\frac{\partial
u^*}{\partial \nu}=0$ on $\partial \Omega$ and we have proved the existence of a non-regular solution solution $u^*$ of
(\ref{10R}) with $\lambda=\lambda^*$.

Let us now show  that for any $\lambda> \lambda^*$ problem
(\ref{10R}) has a non-regular solution solution. Let  $\sigma>1$. Then $\Omega_\sigma:=\{x \in
\mathbb{R}^n:~x\cdot\sigma \in \Omega\} \subset \Omega$, since
$\Omega$ is the star-shaped domain with respect to the origin.
Let us set $ u^*_\sigma(x)=u^{*}(x\cdot\sigma), ~x\in \Omega_\sigma,
$ and  $ u^*_\sigma(x)=0$ in $\Omega\setminus \Omega_\sigma$. Then the following identity
$$
 -\frac{1}{\sigma^2}\,\Delta u_\sigma^* =
 \lambda^* (u_\sigma^*)^{\beta}-(u_\sigma^*)^{\alpha}
$$
weakly holds in $\Omega_\sigma$. Furthermore, since
$u^*_\sigma=0$ and $\frac{\partial u^*_\sigma}{\partial \nu}=0$ on
$\partial \Omega_\sigma$,  this identity weakly holds also in
$\Omega$. This implies that the function
$w(x)=\sigma^\frac{2}{1-\alpha}\cdot u_\sigma(x)$ weakly satisfies
 problem (\ref{10R}) in $\Omega$ with
$\lambda=\sigma^\frac{2(\beta-\alpha)}{1-\alpha} \cdot\lambda^*$.
Note that   $\lambda>\lambda^*$, since $\sigma>1$. This implies
that for any $\lambda\geq \lambda^*$ problem (\ref{10R}) has a
non-regular solution solution.

Let us prove the second part of the Theorem. Note that by Proposition \ref{Lag} 
to any $\lambda \in Z$ it corresponds a weak nonnegative solution $w_\lambda$. Furthermore,  Pohozaev's identity
implies that  $\int \left|\frac{\partial w_\lambda}{\partial
\nu}\right|^2\,x\cdot \nu\,dx>0$, since $P_{\lambda}(w_\lambda)< 0$. Hence there exists a subset $U \subseteq \partial \Omega$
of positive $(n-1)$-dimensional Lebesgue measure such that $\frac{\partial
w_\lambda}{\partial \nu}(s) \neq 0$ for every $s \in U$.

Since $Z$ is an open set, for any $\varepsilon >0$ we can find $\lambda_0 \in Z$ such that $\lambda_0>\lambda^*$ and $\lambda_0-\lambda^*<\varepsilon$. Consider the solution
$w_{\lambda_0}$. Then for any $\lambda>\lambda_0$ the function $w_{\lambda_0}$ is a sub-solution of (\ref{10R}). 

Let us show that (\ref{10R}) has a super-solution. To this end
consider the solution $e \in C^1(\Omega)$ of the following problem
 \begin{equation}
\label{e} \left\{
\begin{array}{l}
\
-\Delta e = 1,
~ x \in \Omega, \\
%\label{1.1}\\ \\
\hspace*{0.2cm} e |_{\partial \Omega}=0,~ x \in \Omega.
%\label{1.2}
%
\end{array}
\right.
\end{equation}
By the maximum
principle for elliptic equations \cite{giltrud} it follows that $e(x)>0$ on $\Omega$ and $\frac{\partial
e}{\partial \nu}(s) <0$ for every $s \in \partial \Omega$. Denote $||e||_\infty=\sup_\Omega|e(x)|$.
Then there exists a sufficient large number $M(\lambda)$ such that
the following inequality
$$
M-\lambda M^\beta||e||_\infty^\beta>0
$$
holds for any $M>M(\lambda)$. Hence and by (\ref{e}) we have
$$
M=-\Delta(Me(x))\geq \lambda (Me(x))^\beta-(Me(x))^\alpha~\mbox{for all}~ x \in \Omega.
$$
Therefore $\overline{u}_\lambda =M\,e$ for any $M>M(\lambda)$ is a super-solution of (\ref{10R}). Furthermore, if $M >M(\lambda)$ is a sufficiently large number, then $\overline{u}_\lambda(x) > w_{\lambda_0}(x)$ in $\Omega$. Thus we may appeal to the method of sub- and super-solutions and therefore there exists a weak solution $w_\lambda \in C^{1,\kappa}(\overline{\Omega})$ for $\kappa \in (0,1)$. The inequality $w_\lambda\geq w_{\lambda_0}$ yields that this solution is not of non-regular solution. Since $\varepsilon>0$ has been taken arbitrary, this completes the proof of the Theorem.
\hspace*{\fill}\rule{3mm}{3mm}\\
\par\noindent

By Corollary \ref{Lag2} we know that $P_{\Lambda_0}(u_{\Lambda_0})<0$. This and Lemma \ref{Pr2} yield that
\begin{cor}
    $\lambda^*<\Lambda_0$.
\end{cor}

\section{Appendix A}
\begin{lem}\label{app}
Assume $\lambda \in [\Lambda_1, +\infty)$ and $u_{\lambda_m}$ is a sequence of solutions of (\ref{min1}),
where $\lambda_m \to \lambda$ as $m\to +\infty$. Then there exists
a subsequence {\rm (}again denoted by $(u_{\lambda_m})${\rm )} and the limit
solution $u_\lambda$ of  (\ref{min1}) such that
$u_{\lambda_m} \to u_\lambda$ strongly
 in $H^1$ as $m \to +\infty$.
\end{lem}
{\it Proof.}\, Let $\lambda \in [\Lambda_1, +\infty)$ and
$u_{\lambda_m}$ be a sequence of solutions of (\ref{min1}), where
$\lambda_m \to \lambda$ as $m\to +\infty$. Let $t_m\geq 0$ and
$v_m\in H^1_0$, $m=1,2,...,$ be such that $u_{\lambda_m}=t_m v_m$,
$||v_m||_1=1$.

As in the proof of Lemma \ref{le1e} using (\ref{Eg1}) it is
derived that $\{t_m\}$ is bounded. This implies that the set
$\{u_{\lambda_m}\}$, $m=1,2,...$ is bounded in $H^1_0$. Hence by
the Sobolev embedding theorem and by the Eberlein-\v{S}mulian
theorem we may assume that $u_{\lambda_{m}} \to \bar{u}_\lambda$
strongly in $L_p(\Omega)$, where $1<p<2^*$, and $u_{\lambda_{m}}
\rightharpoondown \bar{u}_\lambda$ weakly in $H^1_0$ as $ m \to
+\infty$ for some limit point $\bar{u}_\lambda$. This yields that
$\bar{u}_\lambda \in H^1_0$ is a weak nonnegative solution of
(\ref{10R}). As in the proof of Lemma \ref{le1e} using (\ref{Eg1})
it is derived that $\bar{u}_\lambda\neq 0$.

Thus  we have
\begin{eqnarray}
&&E_{\lambda}(\bar{u}_\lambda)\leq \lim_{m\to \infty} E_{\lambda_{m}}(u_{\lambda_{m}}),\label{CFF}\\
&&Q_{\lambda}(\bar{u}_\lambda)\leq 0.
\end{eqnarray}
Let first consider the case $\lambda>\Lambda_1$. By Lemma
\ref{le1e} there exists a solution $u_\lambda$ of (\ref{min1}),
i.e. $u_\lambda \in M_\lambda$ and
$\hat{E}_\lambda=E_\lambda(u_\lambda)$. Then
\begin{eqnarray}
|E_\lambda(u_\lambda)-E_{\lambda_m}(u_\lambda)|<C|\lambda-\lambda_m|,\label{C1}
\end{eqnarray}
where $C<+\infty$. Furthermore, we have
\begin{eqnarray*}
E_{\lambda_m}(u_\lambda)\geq
E_{\lambda_m}(t^2_{\lambda_m}(u_\lambda)u_\lambda)\geq
E_{\lambda_m}(u_{\lambda_m})
 \end{eqnarray*}
provided that $m$ is a sufficiently large number. Thus by (\ref{C1}) we have
\begin{eqnarray*}
E_\lambda(u_\lambda)+C|\lambda-\lambda_m|
>E_{\lambda_m}(u_\lambda)\geq
E_{\lambda_m}(u_{\lambda_m}),
\end{eqnarray*}
and therefore $\hat{E}_\lambda:=E_\lambda(u_\lambda)\geq \lim_{m\to \infty} E_{\lambda_{m}}(u_{\lambda_{m}})$. Using now (\ref{CFF}) we deduce
$$
E_{\lambda}(\bar{u}_\lambda)\leq \hat{E}_\lambda.
$$
Assume $Q_\lambda(\bar{u}_\lambda) < 0$. Then
$Q_\lambda(t^2_\lambda(\bar{u}_\lambda)\bar{u}_\lambda) = 0$, i.e.
$t^2_\lambda(\bar{u}_\lambda)\bar{u}_\lambda\in M_\lambda$ and
$E_\lambda(t^2_\lambda(\bar{u}_\lambda)\bar{u}_\lambda)<E_\lambda(\bar{u}_\lambda)\leq
\hat{E}_\lambda$. Hence we get a contradiction and therefore
$E_\lambda(\bar{u}_\lambda)= \hat{E}_\lambda$, $Q_\lambda(\bar{u}_\lambda) =
0$. Furthermore, $u_{\lambda_m} \to \bar{u}_\lambda$ strongly
 in $H^1$ as $m \to +\infty$ since $Q_\lambda(\bar{u}_\lambda)
= 0$ and hence, we get the proof of the lemma in the case $\lambda >\Lambda_1$.

Assume now that $\lambda =\Lambda_1$.
Then  $\Lambda_1<\lambda_{m}$. By
definition (\ref{P20}) we see that
$\Lambda_1<\lambda_1(u_{\lambda_{m}})\leq \lambda_{m}$ and
therefore $\lambda_1(u_{\lambda_m}) \to \Lambda_1$ as $m \to
\infty$. Thus $(u_{\lambda_m})$ is the minimizing sequence of
(\ref{P30}) and as  above we deduce that $\lambda_1(\bar{u}_{\lambda})= \Lambda_1$.
This implies the proof of the lemma when $\lambda =\Lambda_1$.

\hspace*{\fill}\rule{3mm}{3mm}\\
\par
\end{document}